\documentclass[runningheads]{llncs}
\usepackage{graphicx}
\usepackage{dirtytalk}
\usepackage{babel}
\usepackage[utf8]{inputenc}
\usepackage[T1]{fontenc}
\usepackage{lmodern}
\usepackage{xspace}
\usepackage{mathtools}
\usepackage{amsmath}
\usepackage{enumitem}
\usepackage{accents}
\usepackage{microtype}
\usepackage{hyperref}
\usepackage{cleveref}
\usepackage{csquotes}
\usepackage{pgfplots}
\usepackage{algorithm}
\usepackage{algpseudocode}

\usepackage{amssymb}
\usepackage{hyperref}
\usepackage{todonotes}
\usepackage{color}
\usepackage{algpseudocode}
\usepackage{algorithm}
\usepackage{bbm}
\usepackage{subfig}

\newcommand{\Bcal}{\mathcal{B}}

\newcommand{\Lcal}{\mathcal{L}}
\newcommand{\Ncal}{\mathcal{N}}

\newcommand{\Pcal}{\mathcal{P}}
\newcommand{\Xcal}{\mathcal{X}}
\newcommand{\Ycal}{\mathcal{Y}}
\renewcommand{\P}{\mathbb{P}}
\newcommand{\E}{\mathbb{E}}
\newcommand{\soc}{\text{soc}}

\newcommand{\X}{\mathcal{X}}
\newcommand{\R}{\mathbb{R}}

\newcommand{\abbr}[1][abbrev]{#1.\xspace}

\newcommand{\ie}{\abbr[i.e]}

\algblock{Input}{EndInput}
\algnotext{EndInput}
\algblock{Output}{EndOutput}
\algnotext{EndOutput}

\usepackage{soul}
\usepackage{tikz}
\usetikzlibrary{calc}

\makeatletter
\newif\if@anonymize

\@anonymizefalse  

\if@anonymize
  \newcommand{\highlight@DoHighlight}{
    \fill [outer sep = -15pt, inner sep = 0pt, color=black]
          ($(begin highlight)+(0,8pt)$) rectangle ($(end highlight)+(0,-3pt)$) ;
  }

  \newcommand{\highlight@BeginHighlight}{
    \coordinate (begin highlight) at (0,0) ;
  }

  \newcommand{\highlight@EndHighlight}{
    \coordinate (end highlight) at (0,0) ;
  }

  \newdimen\highlight@previous
  \newdimen\highlight@current
  \newlength{\item@width}

  \DeclareRobustCommand*\anonymize{%
    \SOUL@setup
    \def\SOUL@preamble{%
      \begin{tikzpicture}[overlay, remember picture]
        \highlight@BeginHighlight
        \highlight@EndHighlight
      \end{tikzpicture}%
    }%
    \def\SOUL@postamble{%
      \begin{tikzpicture}[overlay, remember picture]
        \highlight@EndHighlight
        \highlight@DoHighlight
      \end{tikzpicture}%
    }%
    \def\SOUL@everyhyphen{%
      \discretionary{%
        \SOUL@setkern\SOUL@hyphkern
        \SOUL@sethyphenchar
        \tikz[overlay, remember picture] \highlight@EndHighlight ;%
      }{%
      }{%
        \SOUL@setkern\SOUL@charkern
      }%
    }%
    \def\SOUL@everyexhyphen##1{%
      \SOUL@setkern\SOUL@hyphkern
      \settowidth{\item@width}{##1}%
      \makebox[\item@width]{}%
      \discretionary{%
        \tikz[overlay, remember picture] \highlight@EndHighlight ;%
      }{%
      }{%
        \SOUL@setkern\SOUL@charkern
      }%
    }%
    \def\SOUL@everysyllable{%
      \begin{tikzpicture}[overlay, remember picture]
        \path let \p0 = (begin highlight), \p1 = (0,0) in \pgfextra
          \global\highlight@previous=\y0
          \global\highlight@current =\y1
        \endpgfextra (0,0) ;
        \ifdim\highlight@current < \highlight@previous
          \highlight@DoHighlight
          \highlight@BeginHighlight
        \fi
      \end{tikzpicture}%
      \settowidth{\item@width}{\the\SOUL@syllable}%
      \makebox[\item@width]{}%
      \tikz[overlay, remember picture] \highlight@EndHighlight ;%
    }%
    \SOUL@
  }
\else
  \newcommand{\anonymize}[1]{#1}
\fi

\makeatother

\newcommand{\review}{\textcolor{black}} 

\begin{document}
\title{Alternating mixed-integer programming and neural network training for approximating stochastic two-stage problems\thanks{\anonymize{Supported by  Digital Futures and C3.ai Digital Transformation Institute.}}}

 \titlerunning{Alternating MIP and NN training for approximating 2SP} 
%
\author{\anonymize{Jan Kronqvist}\inst{1}\orcidID{\anonymize{0000}-\anonymize{0003}-\anonymize{0299}-\anonymize{5745}} \and \anonymize{Boda Li}\inst{2}\orcidID{\anonymize{0000}-\anonymize{0002}-\anonymize{7934}-\anonymize{9920}} \and \anonymize{Jan Rolfes}\inst{1}\orcidID{\anonymize{0000}-\anonymize{0002}-\anonymize{5415}-\anonymize{1715}} \and
\anonymize{Shudian Zhao}\inst{1}\orcidID{\anonymize{0000}-\anonymize{0001}-\anonymize{6352}-\anonymize{0968}} }
\authorrunning{\anonymize{Kronqvist et al.}}
%
\institute{\anonymize{KTH - Royal Institute of Technology, Stockholm,
  Sweden, Department of Mathematics, Lindtstedtsvägen 25, SE-100 44 Stockholm;
  Sweden}\\
\email{\anonymize{\{jankr,jrolfes,shudian\}@kth.se}} \and 
\anonymize{ABB Corporate Research Center, Wallstadter Str. 59, 68526 Ladenburg, Germany \\ \email{\{davidlee95th\}@gmail.com}}}
\maketitle              
\begin{abstract}
The presented work addresses two-stage stochastic programs (2SPs), a broadly applicable model to capture optimization problems subject to uncertain parameters with adjustable decision variables. In case the adjustable or second-stage variables contain discrete decisions, the corresponding 2SPs are known to be $\mathrm{NP}$-complete.
The standard approach of forming a single-stage deterministic equivalent problem can be computationally challenging  even for small instances, as the number of variables and constraints scales with the number of scenarios. To avoid forming a potentially huge MILP problem, we build upon an approach of approximating the expected value of the second-stage problem by a neural network (NN) and encoding the resulting NN into the first-stage problem. The proposed algorithm alternates between optimizing the first-stage variables and retraining the NN. We demonstrate the value of our approach with the example of computing operating points in power systems by showing that the alternating approach provides improved first-stage decisions and a tighter approximation between the expected objective and its neural network approximation.
\end{abstract}

\section{Introduction}
The area of mathematical optimization often assumes perfect information on the parameters defining the optimization problem. However, in many real-world applications, the input data may be uncertain. The uncertain parameters when optimizing power systems can, e.g., be the power output of wind and solar power units or hourly electricity demands. To address this challenge, two major approaches have been developed in the current literature, namely, robust optimization (RO), see, e.g., \cite{bertsimas2011theory} for further details and stochastic optimization (SO), see, e.g., the seminal surveys \cite{Birge2011a} and \cite{Grossmann2021a}. In RO, only limited information on the underlying distribution of uncertainties is available, e.g., that the uncertain parameter vector belongs to a polyhedron, and the goal is to protect against the worst possible outcome of the uncertainty. In SO, it is typically assumed that the underlying probability distribution of the uncertain parameters is known, or at least that representative samples are available. Here, often the goal is to optimize over the expectation of the uncertain parameters, which results in less conservative solutions compared to RO. The ability to handle uncertainties makes both methods valuable in power-system-related optimization. A comparison between RO and SO for unit commitment with uncertain renewable production is presented in \cite{kazemzadeh2019robust}. However, SO tends to be more prevalent in power system operation optimization, as operating scenarios are readily available. Applications of SO include the optimal operation of high-proportion renewable energy systems \cite{meibom2010stochastic,sharafi2015stochastic,yu2019stochastic}, pricing of energy and reserves in electricity markets \cite{4162624,chazarra2016stochastic,zhang2009stochastic} or day-ahead power scheduling \cite{6714594,shams2018stochastic,mohseni2019voltage} -- to name a few.

Moreover, in some real-world applications, not all the decisions have to be taken before the uncertainty realizes. Here, we distinguish between \emph{first-stage} or here-and-now decisions and \emph{second-stage} or wait-and-see decisions. These \emph{stochastic two-stage} problems or \emph{2SP}s have drawn significant attention over the last years, both due to the practical applicability and computational challenges. Methods for dealing with 2SPs include L-shaped methods \cite{laporte1993integer}, Lagrangian relaxation-based methods \cite{caroe1999dual}, and the sample average approach of Kleywegt et al. \cite{Kleywegt2002a}. For more information on 2SP and algorithmic approaches, we refer to \cite{Grossman2022a,kuccukyavuz2017,Dumouchelle2022a} and the references therein. 

Here, we focus on 2SPs, where the uncertainties either belong to a bounded discrete distribution or can be approximated sufficiently close by such a distribution. Consequently, we can restrict ourselves to a finite set of scenarios.
Furthermore, we are assuming some of the second-stage variables are restricted to integer values and that all constraints are linear. For such problems, it is possible to represent the 2SP as a classical single-stage mixed-integer linear program (MILP) through the so-called deterministic equivalent \cite{Birge2011a}. The deterministic equivalent formulation creates one copy of the second-stage variables and constraints for each scenario, thus resulting in a potentially huge MILP problem. Especially if the second stage involves many integer variables, it may be computationally intractable to consider a large number of scenarios. To overcome this issue, Dumouchelle et al. \cite{Dumouchelle2022a} proposed \review{the so-called Neur2SP approach} ,where the value function, i.e., the optimal objective value of the second stage, is approximated by a neural network (NN). By selecting a MILP representable NN architecture, e.g., ReLU activation functions \cite{lomuscio2017approach,fischetti2018deep,tsay2021partition}, the NN approximation of the second-stage problem along with the first-stage problem can then be integrated into a single stage MILP. The advantage of this approach is that the MILP presentation of the NN can require far fewer variables and constraints than forming the deterministic equivalent. The potential drawback is that we only utilize an approximation of the expected value of the second-stage problem. 

We build upon the approach presented by Dumouchelle et al. \cite{Dumouchelle2022a} and investigate the proposed algorithm for an application in optimizing power systems. We show that proper training of the NN is crucial, and we propose a dynamic sampling as a training method. Here, we directly train a NN to approximate the expected value function. Our approach can be viewed as a surrogate-based optimization \cite{bhosekar2018advances,forrester2009recent} approach for 2SPs, as we are iteratively forming a surrogate of the expected value function which is globally optimized to find a new solution and to improve the approximated expected value function. The proposed algorithm is not limited to the specific application in Section~3, but applicable to any 2SP of suitable structure, e.g., see the applications used in \cite{Dumouchelle2022a}. 

\review{For the application of finding optimal operating points of a power grid considered in the present article, the existing approaches in the power system literature usually approximate the integer variables in the second stage, e.g. by ignoring the action variable of energy storage, see \cite{Cao2020a} or only search for a local optimum of the model, e.g. by applying an augmented Benders decomposition algorithm, see \cite{Cao2022a}. Thus, the Neur2SP approach keeps the underlying integral structure of the second-level, while at the same time aims for a global optimum.}

The outline of the present article is as follows. In Section \ref{Sec: problem setting}, we present the stochastic programming formulation \cite{Dumouchelle2022a} and discuss our methodology. Then, we connect this modeling to achieve data-driven operating states of smart converters in power systems in Section \ref{Sec: powerflow_theory}. Subsequently, we present a numerical comparison between these data-driven results with a baseline algorithm based on the Neur2SP framework in Section \ref{Sec:computationals_results}.

\section{Replacing the second-stage in stochastic programming by a neural network}\label{Sec: problem setting}
The two-stage stochastic programming consists of two types or \emph{levels} of decision variables. The \emph{first-level} variables, denoted by $x\in \Xcal$, describe an initial planning approach. Then, a random vector $\xi\in\Omega$ supported on a compact domain $\Omega\subseteq \R^I$ and distributed by a probability distribution $\P\in \Pcal(\Omega)$ affects the outcome of this initial planning. However, one can adjust the initial planning $x$ after the random vector $\xi$ realizes by choosing the \emph{second-level or recourse} variables $y\in \Ycal(x,\xi)$, wherein the present article, we suppose that $\Ycal(x,\xi)\subseteq \R^n\times \{0,1\}^l$ is assumed to be defined through linear constraints. Thus, we can summarize a \emph{stochastic two-stage problem} as follows:
\begin{equation}\label{Prob: entire_adjustable_SO}
    \min_{x\in \Xcal}~G(x) + \E_\P \left(\min_{y\in \Ycal(x,\xi)} c^\top y\right),
\end{equation}
where in this article $\Xcal$ and $\Omega$ are assumed to be polytopes, \review{$c$ denotes the objective coefficients}, and $\mathcal{Y}(x,\xi)$ a polytope intersected with an integer lattice.
Note, that \eqref{Prob: entire_adjustable_SO} is considered to be a challenging problem as it contains the $\mathrm{NP}$-complete MIP $\min_{y\in \Ycal(x,\xi)} c^\top y$ as a subproblem. Moreover, stochastic two-stage problems are in general notoriously challenging, see e.g., Section 3 in \cite{Grossmann2021a} for an overview.

In the present article, we restrict ourselves to the estimated expected value
\begin{equation}\label{Prob: entire_adjustable_SO_sampled}
    \min_{x\in \Xcal}~G(x) + \frac{1}{m}\sum_{j=1}^m \left(\min_{y\in \Ycal(x,\xi_j)} c^\top y\right),
\end{equation}
where $\xi_j$ denote scenarios sampled from $\P$. For the numerical experiments in Section \ref{Sec:computationals_results}, we limit the number of scenarios to $m=400$. We then aim to learn the mapping $$Q(x)\coloneqq \frac{1}{m}\sum_{j=1}^m\left(\min_{y\in \Ycal(x,\xi_j)} c^\top y\right)$$ 
by training a series of neural networks inspired by the Neur2SP framework developed in \cite{Dumouchelle2022a}. However, in contrast to the approach in \cite{Dumouchelle2022a}, we aim to learn the mapping $x\mapsto Q(x)$, by training a series of neural networks with iteratively generated training data instead of predefined data, \review{see Algorithm~\ref{alg:our_algorithm}}.


\subsection{Embedding the neural network into 2SP}\label{subsec:nn_encode}
We approximate the map $x\rightarrow Q(x)$ by training a fully-connected neural network with 2-hidden layers, where both layers have 40 neurons, denoted as NN $2\times 40$. \review{Hence, the relation between $l$-th layer input  $x^l \in \mathbb{R}^{n_l}$ and output $x^{l+1} \in \mathbb{R}^{n_{l+1}}$ in such neural network with the ReLu activation function is} 
\begin{equation}\label{eq:relu}
x^{l+1} = \max \{0, W^lx^l + b^l\},
\end{equation}
where $W^l \in \R^{n_{l+1}\times n_{l}}$ is the weight matrix and $b^l \in \R^{n_{l+1}}$ is the bias vector.
The ReLu activation function is piece-wise linear and the big-M formulation was the first presented to encode it and is still common~\cite{fischetti2018deep,lomuscio2017approach}. In this way, with binary variable $\sigma^l \in \mathbb{R}^{n_{l+1}}$, \eqref{eq:relu} is equivalent to

\begin{equation}
\begin{aligned}
\label{eq:big-M}
    (w_i^l)^\top x^l + b^l_i \leq x^{l+1}_i, \\
   (w^l_i)^\top x^l + b^l_i - (1-\sigma_i^l)LB^{l+1}_i \geq x^{l+1}_i,\\
    x^{l+1}_i \leq \sigma^l_i UB^{l+1}_i,\\
    \sigma^l_i \in \{0,1\},~x^{l+1}_i\geq 0, \forall i \in [n_{l+1}],
\end{aligned}
\end{equation}
where $w^l_i$ is the $i$-th row vector of $W^l$, $b^l_i$ is the $i$-th entry of $b^l$,  $LB^{l+1}_i$ and $UB^{l+1}_i$ are upper and lower bounds on the pre-activation function over $x^{l+1}_i$, such that $LB^{l+1}_i \leq (w_i^l)^\top x^l + b^l_i \leq UB^{l+1}_i$. Encoding ReLu functions has been an active research topic in recent years, as it enables a wide range of applications of mixed-integer programming to analyze and enhance NNs, such as generating adversarial examples \cite{fischetti2018deep}, selecting optimal inputs of the training data~\cite{kronqvist2022p,zhao2023model}, lossless compression \cite{serra2020lossless}, and robust training~\cite{raghunathan2018certified}.  

By encoding the trained NN as \eqref{eq:big-M}, we can approximate \eqref{Prob: entire_adjustable_SO} by the following MIP:

\begin{subequations}\label{Prob: 2SP_with_NN}
\begin{align}
    \min~& G(x^1) + x^L \\
    \text{s.t.}~& W^l x^l + b^l  \leq x^{l+1}, &&\forall~l \in [L-1], \label{eq:relu_1} \\
     &W^l x^l + b^l - \textrm{diag}(LB^{l+1})(\textbf{1}-\sigma^{l+1}) \geq x^{l+1}, &&\forall~l \in [L-1],\label{eq:relu_2}\\
     & x^l\leq \textrm{diag}(UB^l)\sigma^l , &&\forall 
     ~l \in \{2,\dots ,L\}, \\
     & x^{L} = W^{L-1}x^{L-1}+b^{L-1},\\
    &\sigma^l \in \{0,1\}^{n_l},&& \forall 
     ~l \in \{2,\dots ,L\},\label{eq:relu_3}\\
     & x^l \in \mathbb{R}^{n_l}_+,&&\forall~l \in [L-1],\label{eq:relu_var_size} \\
    & x^1\in \Xcal,x^{L}\in \R,
\end{align}
\end{subequations}
where \review{$x^1\coloneqq x$} is the input-layer variable that belongs to the polytope $\Xcal$, and $x^L$ is the output of the neural network which approximates the expected value $Q(x^1)$.

\subsection{Alternating MIP and NN training}

Since the approximation quality of \eqref{Prob: 2SP_with_NN} with regards to \eqref{Prob: entire_adjustable_SO_sampled} highly depends on the accuracy of the neural network, the approach has significant drawbacks for problems with high dimensional first-level decisions, i.e., high-dimensional $\X$. This drawback becomes even more significant as optimal solutions for \eqref{Prob: entire_adjustable_SO_sampled} are often attained at the boundary or even vertices of $\X$. Here, the approximation by the neural network tends to be worse due to a lack of training data around the vertices. Note that, sampling around each vertex is typically computationally intractable as there are potentially exponentially many such vertices.

The approach we propose in the present article uses intermediate solutions of \eqref{Prob: 2SP_with_NN} to inform the retraining of neural networks. To this end, after training the NN with an initial uniformly distributed sample of $\X$, we solve \eqref{Prob: 2SP_with_NN} and add sample points around the computed optimal solution $x^*_k$ to our training data. Algorithm~\ref{alg:our_algorithm} illustrates our approach.

\begin{algorithm}
\caption{Alternating MIP-NN Algorithm for 2SP (MIP-NN 2SP)}\label{alg:our_algorithm}
\begin{algorithmic}[1]
\Input\ Initial training data $\bar{\X}$, $\alpha=0.99$\;
\EndInput
\Output\ Approximately optimal solution $x^*_{100}$ for \eqref{Prob: entire_adjustable_SO}
\EndOutput

\State Train a NN $2\times 40$ with training data set $\Bar{\X}$\;
\For{$k = 1,\ldots , 100$}
    \State $x^*_k\gets \text{argmin}$ \eqref{Prob: 2SP_with_NN} \label{line:mip} \;
    \For{$i =1,\dots,50$} \label{line:boundary_gen_loop_start}
    \State Sample $x_i \in X$ uniformly\; 
    \State $x_i \gets \alpha x^*_k + (1-\alpha) x_i$ \;
    \State $\bar{\X} \gets \bar{\X} \cup \{x_i\}$ \;
    \EndFor \label{line:boundary_gen_loop_end}
    \State Retrain NN with $\bar{\X}$\;
\EndFor
\end{algorithmic}
\end{algorithm}
We would like to stress, that the number of additional datapoints in every iteration $k$ as well as the parameter $\alpha$ towards $x^*_k$ are chosen arbitrarily. Thus, more careful choices of these parameters may give room for significant improvement of the performance of Algorithm \ref{alg:our_algorithm}.

Here, we would like to \review{briefly comment on the relative computational expenses of the different parts of Algorithm~\ref{alg:our_algorithm}. We specify the considered application as well as our computational ressources in Section \ref{Sec:computationals_results}.} The most time-consuming computation of Algorithm~\ref{alg:our_algorithm} includes solving optimization problems \eqref{Prob: 2SP_with_NN} and $\min_{y\in \Ycal(x,\xi)} c^\top y$. 
With the default settings, Gurobi \review{\cite{gurobi}} can solve \eqref{Prob: 2SP_with_NN} encoding an NN $2\times 40$ within 5 seconds \review{on a standard notebook}. Moreover, the generation of new training data (see line~\ref{line:boundary_gen_loop_start}-\ref{line:boundary_gen_loop_end}), i.e., the optimization for each scenario $\xi$, can be solved in parallel, where solving the problem~$\min_{y\in \Ycal(x,\xi)} c^\top y$ with a pair $x$ and $\xi$ only takes 0.03 seconds. 

\section{Application to smart converters in power system networks}\label{Sec: powerflow_theory}

In order to assess the practical value of the above approach, we apply it to a power system with smart inverters and energy storage. Its effects include reducing power generation costs and improving operational performance are then observed. Since the proposed heuristic given by Algorithm \ref{alg:our_algorithm} may compute first-stage decisions $x$ that lead to an empty second-stage, i.e., $\Ycal (x,\xi)=\emptyset$, we assume in the remainder of this article that every feasible $x\in \X$ leads to non-empty $\Ycal(x,\xi)$ for every scenario $\xi$. \review{In particular, this assumption is valid whenever, we discuss a value function $Q:\X \times \Omega \rightarrow \R, Q(x,\xi)\coloneqq \min_{y\in \Ycal(x,\xi)} c^\top y$ as is often the case in the existing literature, see e.g. \cite{Dumouchelle2022a}, and also justified as the discussion in Section 2.1 in \cite{Shapiro2007a} illustrates.}

\review{Nevertheless}, it is crucial to verify whether this assumption actually holds. \review{For the power system application studied below however, numerical experiments indicate, that} a feasible, but potentially costly $y\in \Ycal(x,\xi)$, exists and thus the use of neural networks will not have an adverse impact on the feasibility of the problem. \review{Moreover, from an engineering perspective, the power system has a slack node, see $b_1$ in Figure \ref{Fig: 5-bus_system}, that ensures stable operation when the system deviates from the set points, although again the operating cost may increase significantly}. On the contrary, it can effectively improve the computational efficiency of the solution. The reasons are as follows:

\begin{enumerate}

\item This paper establishes a two-stage optimization model to solve the day-ahead scheduling strategy of the power system (i.e., the first-stage strategy). We only use neural networks to replace the optimization model in the second stage. The first-stage decisions still satisfy the corresponding hard constraints (e.g., power balance).

\item The second-stage optimization model corresponds to the intra-day operation of the power system. In the cases studied by this paper, the system is connected to the external grid and has energy storage inside as a buffer to balance the energy. From the perspective of practical engineering experience, when the decisions of the first stage are within a reasonable range, the operation strategy of the second stage is usually feasible. Therefore, it is reasonable to use the first-stage decision variables to estimate the operating cost of the second stage through a neural network.

\item In practical operations, the day-ahead operation (i.e., the first stage) needs to obtain results within tens of minutes. When there are a large number of integer variables in the second stage of the 2SP (energy storage requires frequent changes in operational states, which introduces a large number of binary variables into the second stage), this runtime requirement is difficult to guarantee. However, since replacing the second stage of \eqref{Prob: entire_adjustable_SO_sampled} by a neural network can significantly reduce the runtime, Algorithm~\ref{alg:our_algorithm} has the potential to be applicable in practice.
\end{enumerate}

A simplified DC power flow model, following Kirchhoff laws, is adopted to describe the energy balance and power flow in the power system. We follow closely the notations by \cite{Yang2019a} and the presentation in \cite{Kronqvist2023a} in order to specify the parameters present in \eqref{Prob: entire_adjustable_SO}.

\bigskip

Consider a power grid equipped with a set of buses $\Bcal$ and a set of lines/branches $\Lcal$. Within the grid, power is generated by a set of conventional (fossil fuel) generators denoted by $\Ncal_G$ or by a set of distributed (renewable) generators denoted by $\Ncal_{DG}$. The \emph{system operator} \review{(SysO)} may decide, whether to
\begin{itemize}
    \item store or release power to/from a set of batteries $\Ncal_S$
    \item purchase or sell power on the day-ahead market or intra-day, this power enters/leaves the power system through a trading node with the main grid
\end{itemize}

 The amount of power traded day-ahead, i.e., on the day-ahead or \emph{first-level market} is denoted by $P_{fl}$ and the corresponding market price by $p_{fl}$. Similarly, the amount of power and the market price traded intra-day is denoted by $P_{sl}$ and $p_{sl}$, respectively. A positive value for $P_{fl},P_{sl}$ is interpreted as a purchase, negative values as a sell $P_{fl},P_{sl}$ of energy.
 
The SO's initial planning, refers to deciding the first-level variables $x=(P_G^\top,P_{fl})^\top\in \R^{\Ncal_G}\times \R$ based on the estimated renewable energy production $P_{DG, forecast}\in \R^{\Ncal_{DG}}$ in order to ensure that a given total demand $\sum_{i\in \Bcal} P_{d_i}$ is met. Hence, 
$$\Xcal=\{x\in \R^n:\ \eqref{Constr: fl1_define_Pfl}\ \&\ \eqref{Constr: fl2_bounds_P_G} \},$$
where
\begin{subequations}
    \begin{align}
    & P_{fl}^t = \sum_{i\in \Bcal} P_{d_i}^t - \sum_{i\in \Ncal_G}P_{G_i}^t - \sum_{i\in \Ncal_{DG}} P_{DG_i, forecast}^t & \review{\forall} t\in T, \label{Constr: fl1_define_Pfl}\\
    & P_{G_i,\text{min}}^t \leq P_{G_i}^t \leq P_{G_i,\text{max}}^t & \review{\forall} i\in \Ncal_G,  t\in T,\label{Constr: fl2_bounds_P_G}
\end{align}
\end{subequations}
with given parameters $P_{G,\min}, P_{G,\max}\in \R^{\Ncal_G \times T}$. Note that the \emph{market-clearing condition} \eqref{Constr: fl1_define_Pfl} ensures the active power balance in the whole system. The objective of the first level \review{in \eqref{Prob: entire_adjustable_SO}} is then given by 
$$G(x)=x^\top \text{Diag}(c_2) x + c_1^\top x + c_0,$$
where $c_2,c_1\in \R^{\Ncal_G},c_0\in \R$ are given generator cost parameters. Since the uncertainties will impact the initial planning and may jeopardize the power balance, we consider the capacity of the renewable generators as uncertain, i.e., we denote the second-level variable by $\xi=P_{DG, \text{max}}\in \R^{\Ncal_{DG} \times T}$. As these generators are dependent on weather conditions, which are highly uncertain, this is one of the most common uncertainties faced by modern power grids with a high proportion of renewable energies \cite{pfenninger2014energy,impram2020challenges}. To this end, we draw samples from the following domain
\begin{equation}\label{Eq: Def_Omega}
    \Omega=\left\{P_{DG,\text{max}}\in \R^{\Ncal_{DG}\times T}:\ 0\leq  P_{DG_i, \text{max}}^t \leq P_i^+\ \forall\ i\in \Ncal_{DG}, t\in T \right\},
\end{equation}
where $P_i^+$ denotes the technical limit of the renewable generator, i.e. its capacity under optimal conditions. In particular, since the capacity of the renewable generators is further limited by the weather conditions\review{.}

For our computational results in Section \ref{Sec:computationals_results}, we assume that without these limitations $P_{DG_i,\text{max}}$ are independently distributed according to a normal distribution, i.e. $P_{DG_i,\text{max-wol}} \sim \Ncal(P_{DG_i,forecast}, 0.1\cdot P_{DG_i,forecast})$. Consequently, every realization of $P_{DG_i,\text{max-wol}}$ leads naturally to a sample point 
$$\xi_{DG_i}=P_{DG_i,\text{max}} = \begin{cases} 0 & \text{ if } P_{DG_i,\text{max-wol}} \leq 0\\
 P_{DG_i,\text{max-wol}} & \text{ if } P_{DG_i,\text{max-wol}} \in (0,P_i^+)\\
P_i^+ & \text{ if } P_{DG_i,\text{max-wol}} \geq P_i^+.
\end{cases}$$

On the third level, the SysO adjusts the initial planning according to this realization. In particular, the SysO might regulate the energy output of the conventional generators $P_G$ to $P_{G,\text{reg}}$ by either increasing the production by adding $P_G^+\geq 0$ at a cost $r^+$ or decreasing the production by adding $P_G^-\leq 0$ at a cost $r^-$. Similarly, $P_{DG}\in [0,\xi_{DG}]$ is the adjusted energy production that deviates from its forecast by $P_{DG}^+$ or $P_{DG}^-$ with deviations penalized by $f^+,f^-$ respectively. Moreover, the SysO might also trade power intra-day ($P_{sl}$) or decide to use the batteries by setting ($\mu_{ch},\mu_{dch}\in \{0,1\}^{\Ncal_S}$) and change the charging/discharging quantity $P_{ch},P_{dch}$ in order to balance a potential power deficiency or surplus. Based on these decisions, the following variables vary accordingly:
The \emph{state of charge} of a storage, denoted by $\soc$, the power on a line $(k,l)\in \Lcal$, denoted by $p_{kl}$ and the phase angles of the system, denoted by $\theta$. 
Thus, the SysO adjusts the vector $y=(P_{G_i,\text{reg}}, P_{G_i}^+, P_{G_i}^-, P_{DG_i}, P_{DG_i}^+, P_{DG_i}^-, P_{sl}, P_{ch_i}, P_{dch_i},p_{kl}, \theta, \soc, \mu_{ch}, \mu_{dch} )^\top$ in order to satisfy 
$$y\in \Ycal(x,h)\coloneqq \left\{y\in \R^m:\ \eqref{Constr: adjust_fossil_fuel_generators}-\eqref{Constr: mu_bound}\right\},$$
where the constraints \eqref{Constr: adjust_fossil_fuel_generators}-\eqref{Constr: mu_bound} are given below:

\begin{enumerate}
    \item First, we consider the \emph{Generator and DG output constraints}:
    \begin{subequations}
        \begin{align}
            & P_{G_i,\text{reg}}^t = P_{G_i}^t + P_{G_i}^{t,+}+ P_{G_i}^{t,-} & \forall i\in \Ncal_G, t\in T, \label{Constr: adjust_fossil_fuel_generators}\\
            & P_{G_i,\text{min}}^t \leq P_{G_i,\text{reg}}^t \leq P_{G_i,\text{max}}^t & \forall i\in \Ncal_G, t\in T,\\
            & P_{DG_i}^t = P_{DG_i,\text{forecast}}^t + P_{DG_i}^{t,+}+P_{DG_i}^{t,-} & \forall i\in \Ncal_{DG}, t\in T, \label{Constr: adjust_distributed_generators}\\
            & P_{DG_i,\text{min}}^t \leq P_{DG_i}^t \leq P_{DG_i, \text{max}}^t & \forall i\in \Ncal_{DG}, t\in T, \label{Constr: P_DG_leq_P_DG_max}\\
            & P_{sl}^t-P_{fl}^t=-\mathbbm{1}^\top (P_{G}^{t,+} +P_G^{t,-}) -\mathbbm{1}^\top (P_{DG}^{t,+}+P_{DG}^{t,-})\notag\\
            & \qquad\qquad\qquad\qquad\qquad\qquad\quad\ -\mathbbm{1}^\top (P_{dch}^t-P_{ch}^t) & \forall t\in T.\label{Constr: second_level_market_clearing}
        \end{align}
    \end{subequations}
    Here, Constraints \eqref{Constr: adjust_fossil_fuel_generators}-\eqref{Constr: P_DG_leq_P_DG_max} describe the output range of the conventional and renewable generators. In particular, \eqref{Constr: P_DG_leq_P_DG_max} shows that the third-level variables $P_{DG_i}^t$ are restricted by the uncertainties realized in the second level ($P_{DG_i, \text{max}}^{t}$). The market clearing on the intra-day market, i.e., the actual power demand-supply relations, is reflected by Constraint \eqref{Constr: second_level_market_clearing}.
    \item Second, we consider the \emph{operation constraints}:
    \begin{subequations}
        \begin{align}
            & 0 \leq P_{G_i}^{t,+} \leq P_{G_i, \text{max}}^{t,+} && \forall i\in \Ncal_G, t\in T, \label{Constr: pg1}\\
            & P_{G_i, \text{min}}^{t,-} \leq P_{G_i}^{t,-} \leq 0 && \forall i\in \Ncal_G, t\in T,\\
            & P_{DG_i}^{t,+}\geq 0 && \forall i\in \Ncal_{DG}, t\in T,\\
            & P_{DG_i}^{t,-} \leq 0 && \forall i\in \Ncal_{DG}, t\in T, \label{Constr: pg4}
        \end{align}
    \end{subequations} 
    where Constraints \eqref{Constr: pg1} -\eqref{Constr: pg4} limit the real output derivations of conventional and renewable generators.     
    \item Third, we consider the \emph{power flow constraints}. To this end, we apply a DC approximation for a given line reactance ($x_{ij}>0$) and demand in active power ($P_k^{d,t}$) at every time step $t$ and bus $k$:
    \begin{subequations}
    \begin{align}
        \sum_{l \in \delta(k)} p_{kl}^t =& \sum_{i\in \Ncal_G: i\sim k} P_{G_i,\text{reg}}^t +\sum_{i\in \Ncal_{DG}: i\sim k} P_{DG_i}^t\notag\\
         \qquad\qquad &+ \sum_{i\in \Ncal_S: i\sim k} (P_{dch_i}^t-P_{ch_i}^t) - P_k^{d,t} && \forall k\in \Bcal\setminus\{0\}, t\in T, \label{Constr: power_flow1}\\
         \sum_{l \in \delta(k)} p_{0l}^t =& \sum_{i\in \Ncal_G: i\sim 0} P_{G_i,\text{reg}}^t +\sum_{i\in \Ncal_{DG}: i\sim 0} P_{DG_i}^t \notag\\
         \qquad\qquad &+ \sum_{i\in \Ncal_S: i\sim 0} (P_{dch_i}^t-P_{ch_i}^t) +P_{sl}^t - P_0^{d,t} && \forall t\in T, \label{Constr: power_flow2}\\
         \qquad\qquad p_{ij}^t =& \frac{1}{x_{ij}} (\theta_i^t-\theta_j^t)  &&\forall \{i,j\}\in \Lcal, t\in T, \label{Constr:branch_power_balance}
    \end{align}
    \end{subequations}
    where \eqref{Constr: power_flow1} establishes the nodal power flow balance for every node except the root node, which is addressed separately in \eqref{Constr: power_flow2}. The branch power flow is modeled through \eqref{Constr:branch_power_balance}.     
    \item Fourth, we consider the \emph{branch thermal constraints} that guarantee, that the power flow does not exceed the branch's capacities:
    \begin{subequations}
        \begin{align}
            & p_{ij}^t \leq s_{ij, \max} && \forall \{i,j\}\in \Lcal, t\in T.
        \end{align}
    \end{subequations}    
    \item Fifth, we consider the \emph{storage constraints}. To this end, note that the storage operation involves two actions, the storage action is represented by two binary variables $\mu_{ch}, \mu_{dch}\in \{0,1\}^{\Ncal_S\times T}$, whereas the quantity of the charging/discharging is represented by continuous variables $P_{ch_i}^t, P_{dch_i}^t$. Thus, we have the following constraints:
    \begin{subequations}
        \begin{align}
        & \soc_{i,\text{min}}^t \leq \soc_i^t \leq \soc_{i,\text{max}}^t & \review{\forall} i\in \Ncal_S, t\in T, \label{Constr:soc1}\\
        & \soc_i^t = \soc_i^{t-1} + \frac{(P_{ch_i}^t-P_{dch_i}^t)}{E_i} \Delta T & \review{\forall} i\in \Ncal_S, t\in T,\label{Constr:soc2}\\
        & P_{ch_i}^t, P_{dch_i}^t \geq 0 & \review{\forall} i\in \Ncal_S, t\in T, \label{Constr: P_ch_P_dch_nonneg}\\
        & \mu_{ch_i}^t, \mu_{dch_i}^t \in \{0,1\} & \review{\forall} i\in \Ncal_S, t\in T,\\
        & \mu_{ch_i}^t P_{ch_i,\text{min}}^t \leq P_{ch_i}^t \leq \mu_{ch_i}^t P_{ch_i,\text{max}}^t & \review{\forall} i\in \Ncal_S, t\in T, \label{Constr: P_charge_bound}\\
        & \mu_{dch_i}^t P_{dch_i,\text{min}}^t \leq P_{dch_i}^t \leq \mu_{dch_i}^t P_{dch_i,\text{max}}^t & \review{\forall} i\in \Ncal_S, t\in T,\label{Constr: P_discharge_bound}\\
        & \mu_{ch_i}^t+\mu_{dch_i}^t \leq 1 & \review{\forall} i\in \Ncal_S,  t\in T.\label{Constr: mu_bound}
        \end{align}
    \end{subequations}
\end{enumerate}
Here, the upper/lower bounds of soc, as well as the relationships between soc and charging/discharging actions, are given by Constraints \eqref{Constr:soc1} and \eqref{Constr:soc2}. The connection between the storage actions and the respective quantity is reflected by Constraints \eqref{Constr: P_ch_P_dch_nonneg} -- \eqref{Constr: mu_bound}.

Lastly, the third-level cost function aims to both, minimize the electricity cost as well as reduce the deviation between the intra-day system operation strategy and the day-ahead planning. Thus, the whole adjustment can be summarized as solving
\begin{equation}
    \min_{y\in \Ycal(x,\xi)} c^\top y, \label{eq:first-stage}
\end{equation}
where
$$c^\top y\coloneqq \sum_{t\in T} \sum_{i\in \mathcal{N}_G}(r_i^{+,t}P_{G_i}^{+,t} +r_i^{-,t}P_{G_i}^{-,t}) + p_{sl}^t(P_{sl}^t-P_{fl}^t) + \sum_{i\in\mathcal{N}_{DG}} (f_i^+ P_{DG_i}^{t,+}+f_i^- P_{DG_i}^{t,-}).$$

After having established the model parameters for \eqref{Prob: entire_adjustable_SO}, we continue with a case study in order to demonstrate the practical value of Algorithm \ref{alg:our_algorithm}.

\section{Computational results}\label{Sec:computationals_results}

In this section, we present a case study, where we consider the daily power distribution of a 5-bus instance, based on the ``case5.m'' instance from the matpower library \cite{Zimmerman2011a}. Here, the day is divided into hourly (24 period) time intervals. \review{The MIP problems in Algorithm~\ref{alg:our_algorithm} are solved by Gurobi~\cite{gurobi}, which is the state-of-art MIP solver.}
The computations were executed via Gurobi 10.0.0  under Python 3.7 on a Macbook Pro (2019) notebook with an Intel Core i7 2,8 GHz Quad-core and 16 GB of RAM. The neural network is trained with PyTorch~\cite{NEURIPS2019_9015}.

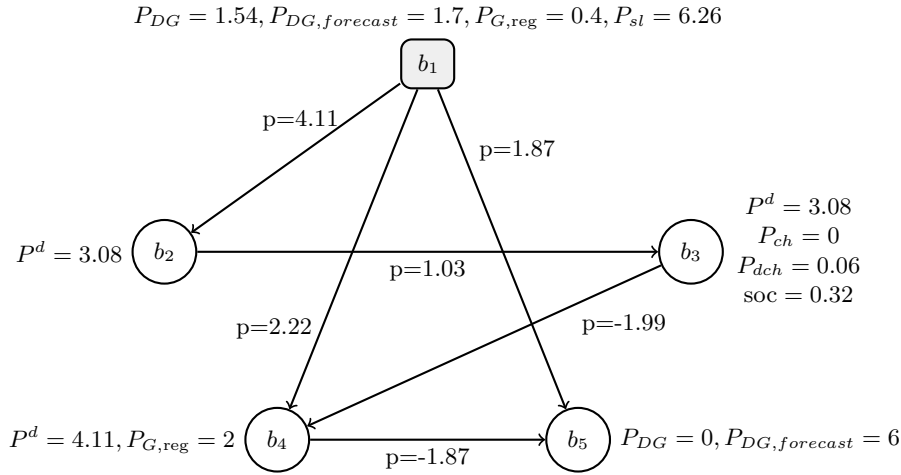
\begin{figure}[ht]

\tikzstyle{block} = [rectangle, draw=black, fill=lightgray!25,text width=1.5em, rounded corners, minimum height=2em, text centered]
		
		\tikzstyle{block3}=[circle, draw=black, fill=white!10,text width=1.5em,rounded corners, minimum height=2em, text centered]
		
		\tikzstyle{line} = [draw, -latex']
		
		\begin{center}
			\begin{tikzpicture}[node distance = 2cm, auto,thick]
				\node [block] at (0,5)(b1) [label={[align=center]above: $P_{DG}=1.54, P_{DG,forecast}=1.7, P_{G,\text{reg}} =0.4, P_{sl}=6.26$}]{$b_1$};
				\node [block3] at (-3.5,2.5) (b2) [label={[align=center]left: $P^d=3.08$}] {$b_2$};
				\node [block3] at (3.5,2.5) (b3) [label={[align=center]right: $\begin{array}{c}
						P^d=3.08 \\
						P_{ch}=0 \\
						P_{dch}=0.06\\
						\text{soc} = 0.32
					\end{array}$}]  {$b_3$};
				\node [block3] at (-2,0) (b4) [label={[align=center]left: $P^d=4.11, P_{G,\text{reg}}=2$}] {$b_4$};
				\node [block3] at (2,0) (b5) [label={[align=center]right: $P_{DG}=0, P_{DG,forecast}=6$}]  {$b_5$};
				\draw [->] (b1) to node [left, near start]{p=4.11} (b2);
				\draw [->] (b1) to node [left, near end]{p=2.22} (b4);
				\draw [->] (b1) to node [auto, near start]{p=1.87} (b5);
				\draw [->] (b2) to node [below, midway]{p=1.03} (b3);
				\draw [->] (b3) to node [auto, near start]{p=-1.99} (b4);
				\draw [->] (b4) to node [below, midway]{p=-1.87} (b5);
			\end{tikzpicture}
			
		\end{center}

\caption[BDD]{The case5.m network with its corresponding generators and an exemplary power flow. $P^d=0$ at buses $1$ and $5$}
\label{Fig: 5-bus_system}
\end{figure}

The network topology of ``case5.m'' is illustrated in Figure \ref{Fig: 5-bus_system}, where we assume bus $1$ to be the slack node, i.e., the trading node connected to other grids. Since the type of generators is not specified in \cite{Zimmerman2011a} and this study solely serves as an academic example, we chose whether a generator in ``case5.m'' is a conventional/renewable one or a storage in the following convenient way: The conventional (fossil fuel) generators are connected to the buses $1$ and $4$, i.e. $\Ncal_G=\{1,4\}$, two distributed generators are connected to buses $1$ and $5$, i.e., $\Ncal_{DG}=\{1,5\}$ and an energy storage unit is connected to bus $3$, i.e. $\Ncal_S=\{3\}$. \review{In order to aid reproducibility and encourage further research on this topic, we provide the underlying system data for public use, see \cite{Kronqvist2023b}.}

Incorporating the instance data closely follows the methodology used in \cite{Kronqvist2023a}. 
Hence, both, the day-ahead and intra-day market prices $p_{fl}, p_{sl}$ were taken as averages from the Pecan street database's \cite{pecanstreet} ``miso'' data set for \review{March 3rd, June 3rd, October 3rd and December 2nd}, 2022. Daily deviations in $P^d$, denoted by $\Delta_d^t$, or daily deviations in $P_{DG,\max}^t$, denoted by $\Delta_{DG}^t$ were taken from the Pecan street database's ``california\_iso'' dataset for \review{the same days}, 2022 in the same way. \review{In particular, we also provide the underlying weather data for public use, see \cite{Kronqvist2023b}.}

Consequently, we model the daily varying demand by $P^d=P^d\cdot \Delta_d^t$ and the varying potential renewable energy production by 
$$P_{DG,\text{forecast}}^t \coloneqq\min\left\{\frac{P_{DG,\min}+P_{DG,\max} }{2}\cdot \Delta_{DG}^t, P_{i,t}^+\right\}. $$

After incorporating this data, we compare the solutions given by Algorithm \ref{alg:our_algorithm} to a baseline experiment, where we replace line \ref{line:boundary_gen_loop_start}-\ref{line:boundary_gen_loop_end} by adding $50$ uniformly distributed sample points drawn from $\X$ in each iteration $k$. This algorithm is inspired by \cite{Dumouchelle2022a} and simply creates the same amount of sample points, but does not incorporate targeted information around $x^*_k$. We chose $3000$ uniformly distributed datapoints from $\X$ as a starting dataset for both algorithms. We use trained neural networks with the same architecture for both experiments, \ie, NN $2\times 40$.


\begin{figure}[!htbp]
	\centering
	\subfloat[Spring]{\includegraphics[width = 2.4in,trim={0cm 0cm 0cm 1.3cm},clip]{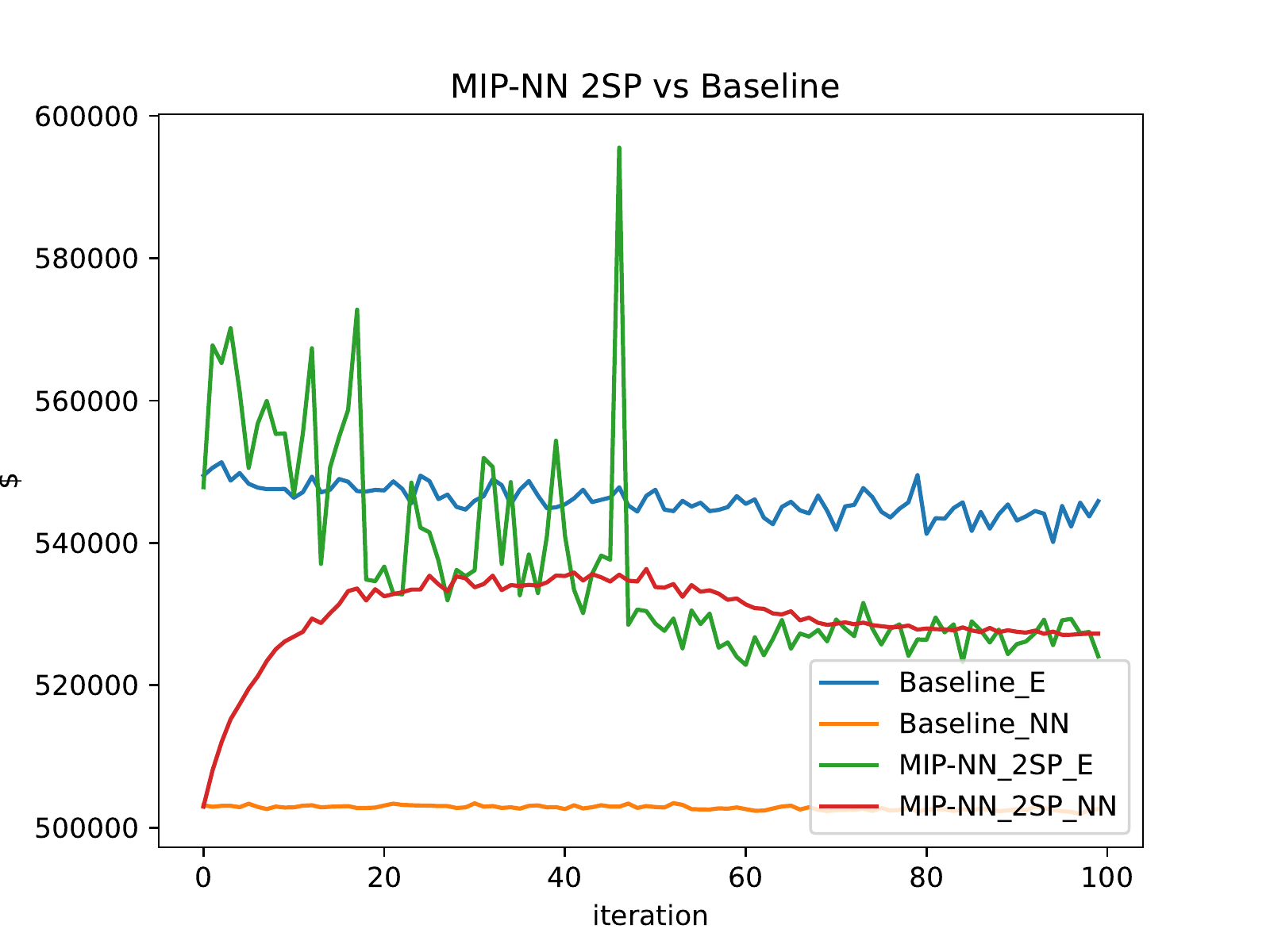}} 
	\subfloat[Summer]{\includegraphics[width = 2.4in,trim={0cm 0cm 0cm 1.3cm},clip]{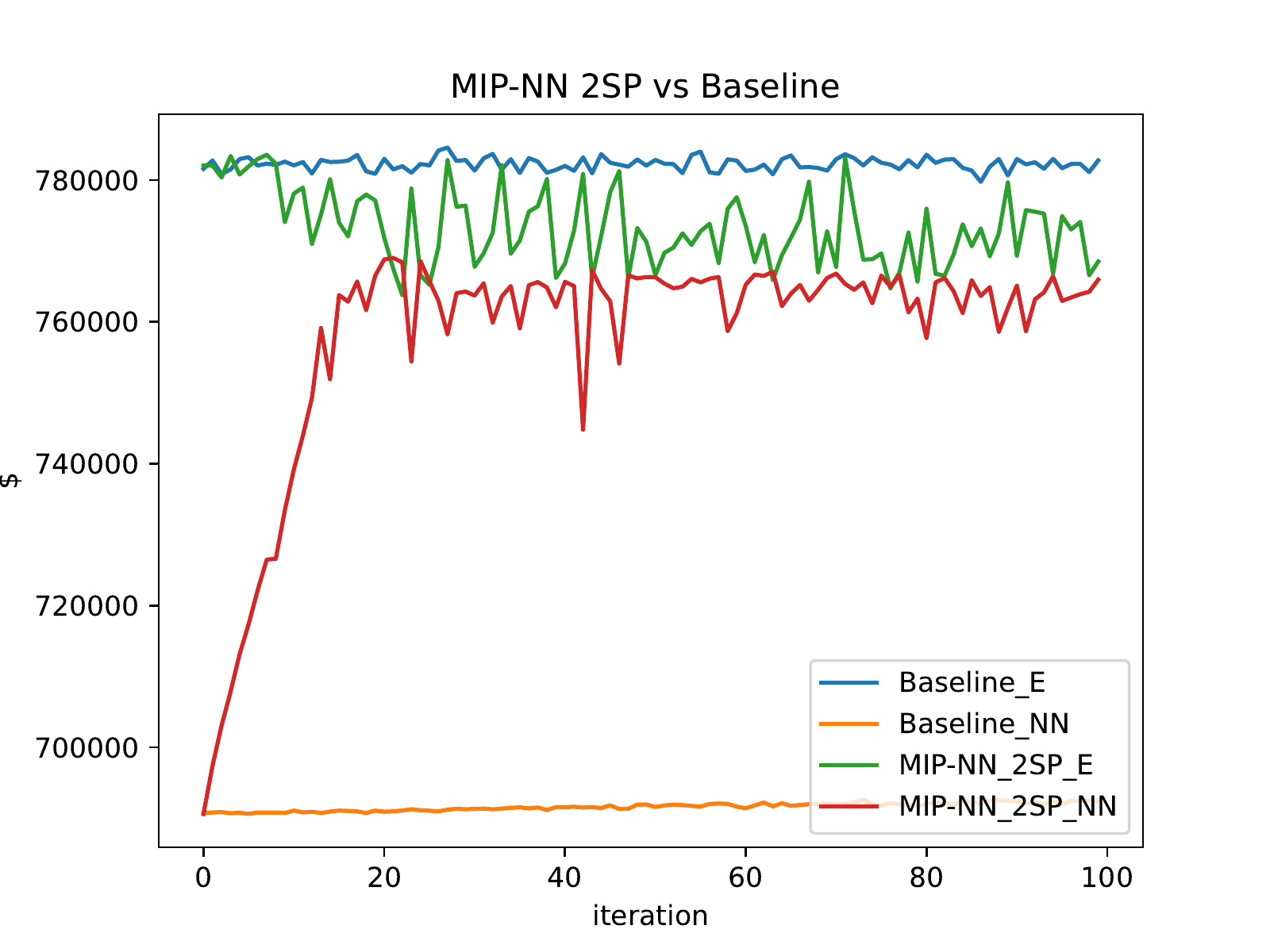}}

 \subfloat[Autumn]{\includegraphics[width = 2.4in,trim={0cm 0cm 0cm 1.3cm},clip]{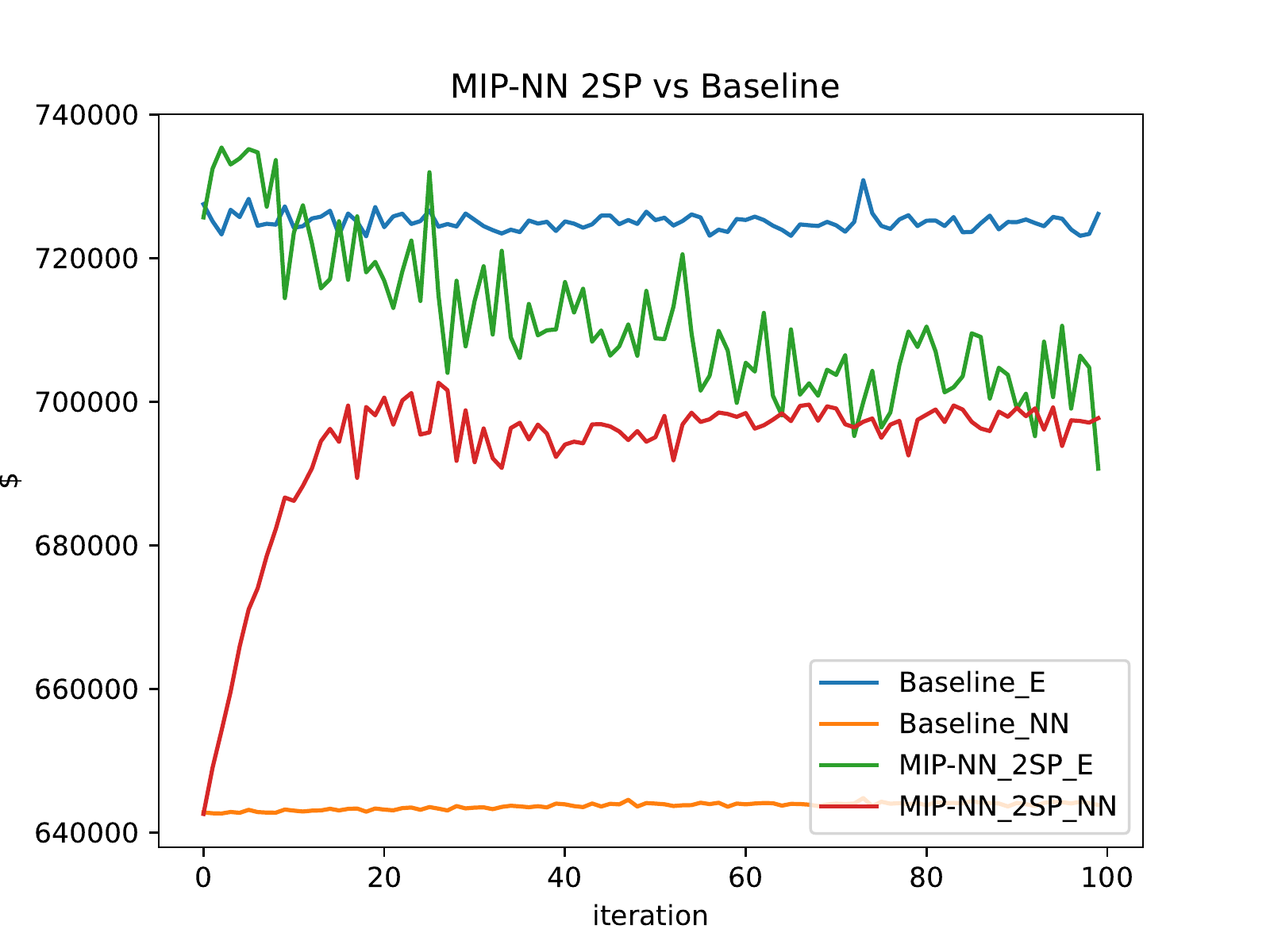}} 
	\subfloat[Winter]{\includegraphics[width = 2.4in,trim={0cm 0cm 0cm 1.3cm},clip]{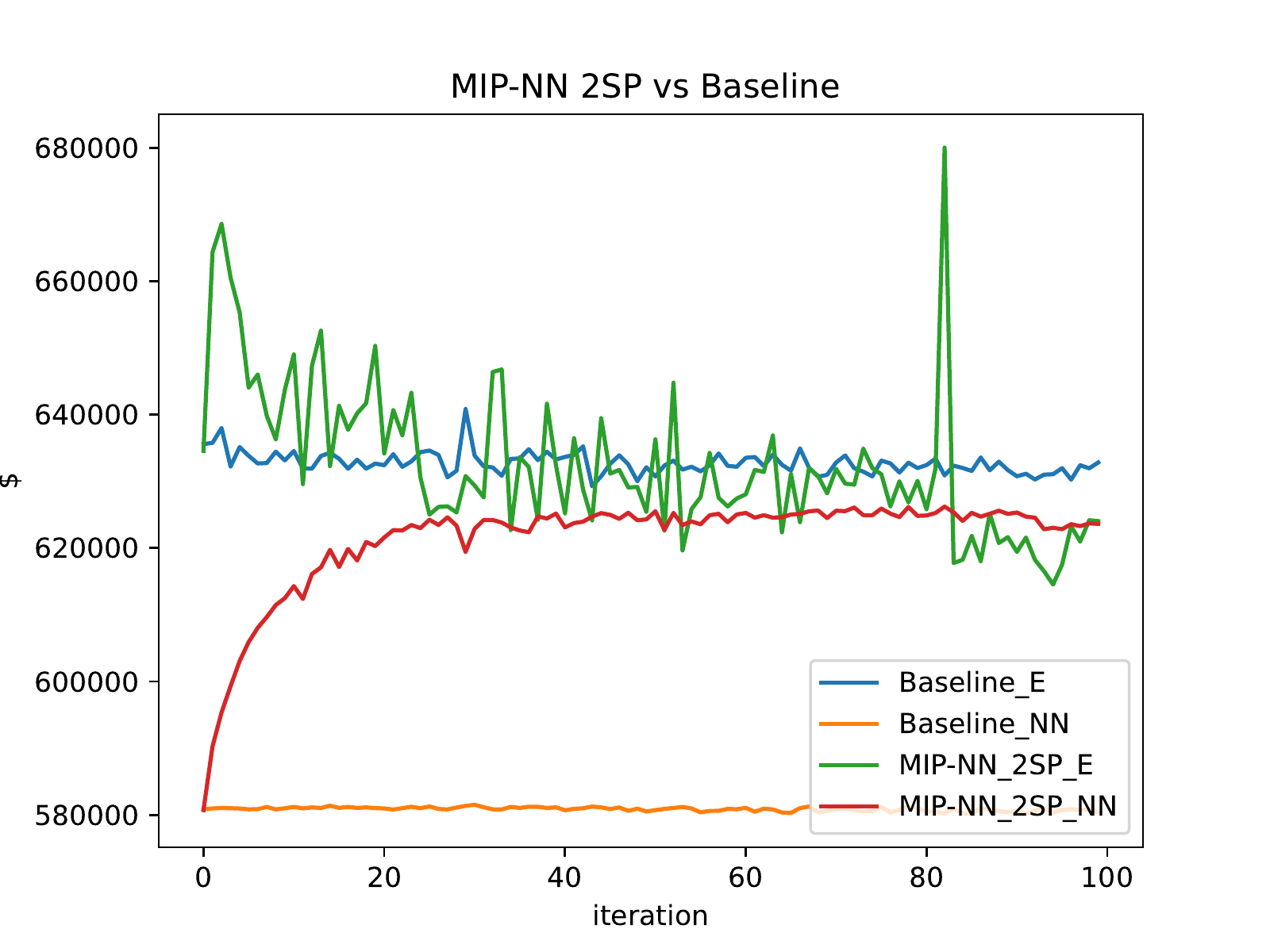}}
	\caption{Comparison of Algorithm \ref{alg:our_algorithm} (in green and red) to successively adding batches of 50 uniformly distributed data points (in blue and orange) with respect to operating costs \review{in \$} for weather data of different seasons.}
	\label{fig:comparison_operating_cost}
\end{figure}

In Figure~\ref{fig:comparison_operating_cost}, MIP-NN\_2SP\_NN (resp. Baseline\_NN) denotes the objective function value of \eqref{Prob: 2SP_with_NN} solved by Algorithm~1 (resp. the baseline experiment) at each iteration. In particular, we include $G(x^*_k)+\frac{1}{400} \sum_{j=1}^{400} Q(x^*_k,\xi_j)$, denoted by MIP-NN\_2SP\_E (resp. Baseline\_E), with optimal solutions $x^*_k$ of line \ref{line:mip} from  Algorithm~1 (resp. the baseline experiment). In this way, the difference between MIP-NN\_2SP\_NN (resp. Baseline\_NN) and MIP-NN\_2SP\_E (resp. Baseline\_E) measures the quality of the approximation by the neural networks at each iteration.

Observe, that after 100 iterations, the baseline experiment cannot close the gap between the estimated expected value of the second stage (\ie, Baseline\_E) and its neural network approximation (\ie, Baseline\_NN) since a gap of $80,000 \$$ remains. In other words, the difference between the NN approximation and the expected value is $80,000 \$$. While including targeted datapoints in the proposed sampling approach significantly reduces the gap between MIP-NN\_2SP\_NN and MIP-NN\_2SP\_E. Thus, showing that the dynamic sampling in Algorithm~\ref{alg:our_algorithm} clearly outperforms sampling uniformly distributed and results in a much more accurate approximation of the expected value. Moreover, the results indicate, that the operating state of the conventional generators computed from the baseline experiment is not optimal, and approximately $20,000 \$$ may be saved by using the operating states computed by Algorithm \ref{alg:our_algorithm}. The results, thus, show the importance of obtaining an accurate approximation of the expected value function. 

\section{Conclusion and Outlook}
We have presented an algorithm for solving 2SP problems with integer second-stage variables by iteratively constructing a NN that approximates the expected value of the second-stage problem. This approach is inspired by the work of \cite{Dumouchelle2022a}. In the algorithm, we propose a dynamic sampling and retraining of the NN to improve the approximation in regions of interest. We numerically evaluate the algorithm in a case study of optimizing a power system. The numerical results highlight the importance of obtaining an accurate approximation of the expected value function and show that a poor approximation can lead to a suboptimal solution. The results strongly support the proposed sampling and retraining approach, and by the smart selection of data points we are able to obtain a much better approximation in the regions of interest compared to a simple uniform sampling. For the baseline, using uniform sampling, the accuracy does not seem to improve even after doubling the number of data points. Whereas, the proposed algorithm quickly obtains an NN that seems to closely approximate the expected value function, at least in the neighborhood of the minimizer of the NN approximated 2SP \eqref{Prob: 2SP_with_NN}. The results, thus, support the idea of approximating the expected value function by an NN and using a MILP encoding of the NN to form a single-stage problem, but clearly show the importance of efficiently training a NN to high accuracy in the regions of interest. 

\review{As Algorithm \ref{alg:our_algorithm} delivers promising numerical results on the considered power system, it is natural to ask, whether these results extend to the benchmark problems given in \cite{Dumouchelle2022a}, i.e. Capacitated Facility location, Investment, Stochastic server location and pooling problems. Moreover, also extensions to ACOPF models similarly to the ones in \cite{Kilwein2021a} may pose interesting challenges for future research.}

\review{\subsubsection{Acknowledgements}
We would like to thank the anonymous referees for their very valuable comments, that helped to significantly improve the quality of this article.} 

%
%
%
\bibliographystyle{splncs04}
\bibliography{mybibliography}

\end{document}